\documentclass[12pt,twoside,french]{amsart}

\usepackage[french]{babel}
\usepackage{amssymb}
\usepackage{graphicx}
\usepackage{amscd}
\usepackage{vmargin}

\setmarginsrb{2.5cm}{2.5cm}{2.7cm}{2cm}{0cm}{1cm}{0cm}{1cm}

\newtheorem{theorem}{Th\'{e}or\`{e}me}
\newtheorem*{theoremeprincipal}{Th\'{e}or\`{e}me principal}

\theoremstyle{plain}

\newtheorem{definition}[theorem]{D\'{e}finition}

\newtheorem{lemma}[theorem]{Lemme}

\newtheorem{proposition}[theorem]{Proposition}

\numberwithin{equation}{section}

\def\preuve{\smallskip\goodbreak{\it Preuve.~~--~\kern.3em}
     \ignorespaces}%
\def\qedbox{$\square$}%
\def\qed{\ifmmode\qedbox\else\unskip\ \hglue0mm\hfill
     \qedbox\smallskip\goodbreak\fi}%
\def\preuvedt#1{\smallskip\goodbreak{\it Preuve du th\'eor\`eme~\ref{#1}.~~--~\kern.3em}
     \ignorespaces}%
\def\idpreuve{\smallskip\goodbreak{\it Id\'ees de la preuve.~~--~\kern.3em}
     \ignorespaces}%\usepackage[all,dvips]{xy}
\def\idpreuvedt#1{\smallskip\goodbreak{\it Id\'ees de la preuve du th\'eor\`eme~\ref{#1}.~~--~\kern.3em}
     \ignorespaces}%
\newcommand{\CC}{\mathbb{C}}
\newcommand{\NN}{\mathbb{N}}
\newcommand{\PP}{\mathbb{P}}

\newcommand{\D}{\mathcal{D}}

\newcommand{\M}{\mathcal{M}}
\newcommand{\I}{\mathcal{I}}
\newcommand{\J}{\mathcal{J}}

\renewcommand{\O}{\mathcal{O}}
\newcommand{\dgl}{$\D$-groupo\"ide de Lie }
\newcommand{\dgls}{$\D$-groupo\"ides de Lie }

\title[Enveloppe galoisienne d'une application rationnelle]{Enveloppe galoisienne \\ d'une application rationnelle de $\mathbb{P}^1$}

\author{Guy Casale}

\address{Laboratoire E. Picard, UMR 5580 UFR MIG, Universit\'e Paul Sabatier \ 31 062 Toulouse cedex 4, France}

\email{casale@picard.ups-tlse.fr}

\begin{document}

\maketitle

%\tableofcontents

\section*{Introduction}

Dans \cite{malg-groupoidedegalois} et \cite{malg-chinois} B. Malgrange d\'efinit la \textit{$\D$-enveloppe}, ou \textit 
{enveloppe galoisienne}, de divers syst\`emes dynamiques analytiques (feuilletage, champ de vecteurs, transformation 
rationnelle \ldots) sur une vari\'et\'e analytique complexe lisse. Heuristiquement, il s'agit du syst\`eme maximal 
d'\'equations aux d\'eriv\'ees partielles portant sur les diff\'eomorphismes locaux de cette vari\'et\'e v\'erifiant les deux 
conditions suivantes :
\begin{itemize}
\item les solutions locales de ce syst\`eme sont stables par inversion et composition,
\item le syst\`eme dynamique est solution (infinit\'esimale) de ce syst\`eme.
\end{itemize} 
Dans le cas des feuilletage, cette objet est fortement li\'e \`a l'existence d'int\'egrales premi\`eres d'un type de 
transcendance particulier (voir \cite{mathese}). Dans le cas des \'equations diff\'erentielles lin\'eaires, la 
$\D$-enveloppe redonne le groupe de Galois diff\'erentiel d'une extension de Picard-Vessiot (\cite{malg-groupoidedegalois}, 
voir aussi \cite{mathese}).\\
Dans cet article nous donnons la liste des transformations rationnelles de $\PP^1$ ayant une 
\textit{$\D$-enveloppe} non triviale, c'est-\`a-dire d\'efinie par au moins une \'equation non nulle. 

\begin{theoremeprincipal}
    \label{theoremeprincipal}
Les seules applications rationnelles de $\PP^1$ dans $\PP^1$ ayant une $\D$-enveloppe non triviale sont, \`a conjugaison 
par une homographie pr\`es, les mon\^omes, les polyn\^omes de Tch\'ebitchev et les exemples de Latt\`es.
\end{theoremeprincipal}

La construction de ces applications sera rappel\'ee dans la premi\`ere partie et les d\'efi\-ni\-tions de $\D$-groupo\"ide de 
Lie et $\D$-enveloppe, dans la deuxi\`eme. 
Nous constatons que cette liste est celle des applications rationnelles admettant un commutant non trivial 
(\cite{RittR,Fatou,Julia,Eremenko}). 
Par analogie avec l'int\'egrabilit\'e des syst\`emes hamiltoniens, les applications rationnelles ayant un commutant non 
trivial sont appel\'ees \textit{``int\'egrables''} (\cite{Veselov}). D'autre part, dans le contexte de cet article l'appellation 
\textit{``int\'egrable''} est justifi\'ee par le fait que plus la \textit{$\D$-enveloppe} est petite moins la dynamique est 
\textit{``transcendante''}.

%%%%%%%%%%%%%%%%%%%%%%%%%%%%%%%%%%%%%%%%%%%%%%%%%%%%%%%%%%%%%%%%%%%%%%%%%
%%%%%%%%%%%%%%%%%%%%%%%%%%%%%%%%%%%%%%%%%%%%%%%%%%%%%%%%%%%%%%%%%%%%%%%%%
%%%%%%                   DEFINITIONS ET ENONCES                     %%%%%
%%%%%%%%%%%%%%%%%%%%%%%%%%%%%%%%%%%%%%%%%%%%%%%%%%%%%%%%%%%%%%%%%%%%%%%%%
%%%%%%%%%%%%%%%%%%%%%%%%%%%%%%%%%%%%%%%%%%%%%%%%%%%%%%%%%%%%%%%%%%%%%%%%%

%%%%%%%%%%%%%%%%%%%%%%%%%%%%%%%%%%%%%%%%%%%%%%%%%%%%%%%%%%%%%%%%%%%%%%%%%
\section{Quelques rappels de dynamique holomorphe \cite{B-M}}
%%%%%%%%%%%%%%%%%%%%%%%%%%%%%%%%%%%%%%%%%%%%%%%%%%%%%%%%%%%%%%%%%%%%%%%%%

%%%%%%%%%%%%%%%%%%%%%%%%%%%%%%%%%%%%%%%%%%%%%%%%%%%%%%%%%%%%%%%%%%%%%%%%%
Soit $R:\PP^1 \to \PP^1$ un application rationnelle de la droite projective. Au voisinage d'un point $p$ fixe ($R(p)=p$), 
r\'epulsif ($|R'(p)|>1$), cette application est lin\'earisable.

\begin{theorem}[de lin\'earisation de K{\oe}nigs]
Soit $f: (\CC,0) \to (\CC,0)$ telle que $|f'(0)|>1$. Il existe une unique application $\Psi : (\CC,0) \to (\CC,0)$ telle 
que $$\Psi'(0)=1 \text{ et }  \Psi^{-1} \circ f \circ \Psi(w)=f'(0)w .$$ On appelle $\Psi$ la lin\'earisante de K{\oe}nigs en $p$.
\end{theorem}

L'ensemble de Julia de $R$ est l'ensemble $\J_R$ des points au voisinage desquels la famille $\{ R^{\circ n}\}_{n \in \NN}$
n'est pas normale. Les points fixes r\'epulsifs des fractions $R^{\circ n}$ appartiennent \`a cet ensemble. Plus 
pr\'ecisement, on a le th\'eor\`eme suivant.

\begin{theorem}
Pour toute fraction rationnelle $R$ de degr\'e au moins deux, les points fixes r\'epulsifs des it\'er\'ees $R^{\circ n}$ 
sont denses dans $\J_R$. 
\end{theorem}

Soit $p$ un point fixe r\'epulsif. La lin\'earisante de K{\oe}nigs s'\'etend par les formules 
$$ R^{\circ n} \circ \Psi (w)=\Psi((R'(p))^n w).$$ L'image de ce prolongement est \'egal \`a $\PP^1 -E$ o\`u $E$ est l'ensemble
exceptionnel, c'est-\`a-dire, l'ensemble de points $q$ dont les pr\'eimages $R^{-n}(\{q\})$ n'accumulent pas l'ensemble de 
Julia de $R$. Lorsqu'il est non vide, cet ensemble est r\'eduit \`a un ou deux points. L'application $\Psi$ d\'efinie un 
rev\`etement ramifi\'e de $\PP^1 -E$ par $\CC$ qui ``semi-conjugue'' $R$ \`a sa partie lin\'eaire en $p$ : 
$$R \circ \Psi = \Psi(R'(p) w).$$

Nous allons maintenant donner des exemples d'applications rationnelles particuli\`eres. La construction de ceux-ci se fait 
\`a partir de leurs lin\'earisantes.\\

\noindent \textbf{Les mon\^omes.} Consid\'erons l'application exponentielle $\exp : \CC \to \PP^1-\{0, \infty\}$ et sur 
$\CC$ la transformation $w \mapsto kw$, o\`u $k$ est un entier. On construit une application rationnelle $M_k$ sur $\PP^1$ 
en posant $M_k(\exp w) = \exp(k w)$. Cette application est \'evidemment l'application monomiale $x^k$ mais son caract\`ere 
rationnel peut se montrer \textit{a priori} en utilisant les formules d'additions de la fonction exponentielle.
 L'ensemble $\J_{M_k}$ est le cercle unit\'e. Le point ``$1$'' est un point fixe r\'epulsif dont la lin\'earisante de 
 K{\oe}nigs se prolonge en l'application exponentielle.

\noindent \textbf{Les polyn\^omes de Tchebitchev.} Consid\'erons cette fois l'application cosinus $\cos : 
\CC \to \PP^1-\{\infty\}$. Des formules d'additions satisfaites par cette application, on peut d\'eduire l'existence de 
polyn\^omes $T_k$ v\'erifiant $T_k(\cos w) = \cos(k w)$. L'ensemble $\J_{T_k}$ est le segment $[ -1,1].$ Le point ``$1$'' est un 
point fixe r\'epulsif dont le lin\'earisante de K{\oe}nigs se prolonge en l'application cosinus.

\noindent \textbf{Les exemples de Latt\`es.} Soit $\wp : \CC \to \PP^1$ la fonction de Weierstrass associ\'ee \`a un r\'eseau
$\Lambda$. Cette fonction satisfait aussi des formules d'additions. Celles-ci permettent de montrer qu'il existe des fractions 
rationnelles $L_k$ satisfaisant $L_k(\wp w) = \wp(k w)$. Plus g\'en\'eralement on appelle exemple de Latt\`es toute 
application rationnelle $L$ telle qu'il existe un tore complexe $\CC/ \Lambda$, une isog\'enie de ce tore $I_\lambda$ et une 
fonction elliptique $p : \CC/\Lambda \to \PP^1$ tel que $p \circ I_\lambda = L\circ p$. L'ensemble $\J_{L}$ est $\PP^1$ car tous 
les points sont ``expansifs''.\\
  
Le th\'eor\`eme suivant d\^u \`a Fatou \cite{Fatou}, Julia \cite{Julia}, Ritt \cite{RittR} et Eremenko \cite{Eremenko} 
(voir aussi \cite{D-S}) donne une caract\'erisation de ces applications rationnelles.
\begin{theorem}
Soient $R_1$ et $R_2$ deux fractions sur $\PP^1$ v\'erifiant 
$$
R_1 \circ R_2 = R_2 \circ R_1
$$
et
$$
R_1^{\circ n_1} \not = R_2^{\circ n_2} \text{ pour tout entier } n_1 \text{ et } n_2.
$$
Alors on est dans un des cas suivants :
\begin{itemize}
\item $R_1$ et $R_2$ sont des mon\^omes \`a multiplication par une racine de l'unit\'e pr\`es,
\item $R_1$ et $R_2$ sont des polyn\^omes de Tch\'ebitchev au signe pr\`es,
\item $R_1$ et $R_2$ sont des exemples de Latt\`es.
\end{itemize} 
\end{theorem}

%%%%%%%%%%%%%%%%%%%%%%%%%%%%%%%%%%%%%%%%%%%%%%%%%%%%%%%%%%%%%%%%%%%%%%%%%
\section{Quelques rappels sur la notion de $\D$-groupo\"ide de Lie sur $\PP^1$}
%%%%%%%%%%%%%%%%%%%%%%%%%%%%%%%%%%%%%%%%%%%%%%%%%%%%%%%%%%%%%%%%%%%%%%%%%

Pour une introduction aux notions de $\D$-groupo\"ide de Lie et $\D$-alg\`ebre de Lie, le lecteur pourra consulter \cite{malg-chinois}
et pour plus de d\'etails \cite{malg-groupoidedegalois}. 
%%%%%%%%%%%%%%%%%%%%%%%%%%%%%%%%%%%%%%%%%%%%%%%%%%%%%%%%%%%%%%%%%%%%%%%%%
\subsection{Structure differentielle et alg\'ebrique de $J^*(\PP^1)$ }
\text{   \\}

Nous noterons $J_k^*(\PP^1)$ l'espace des jets d'ordre $k$ d'applications inversibles de $\PP^1$ dans $\PP^1$. Si on choisit
deux cartes $(U,x)$ et $(V,y)$  de $\PP^1$, un jet d'ordre $k$ s'\'ecrit $(x,y,y_1,\ldots,y_k)$ avec $y_1\not = 0$. 
L'anneau $\O(U \times V)[y_1,y_1^{-1},\ldots,y_k]$ est l'anneau des \'equations diff\'erentielles d'ordre $k$ sur les jets 
d'applications inversibles de $U$ dans $V$. On munit ainsi $\PP^1 \times \PP^1$ d'un faisceau d'anneaux $\O_{J_k^*(\PP^1)}$. 
Les inclusions naturelles $\O_{J_k^*(\PP^1)} \subset \O_{J_{k+1}^*(\PP^1)}$ permettent de d\'efinir $\O_{J^*(\PP^1)} = 
\lim \O_{J_k^*(\PP^1)} $ le faisceau d'anneaux des \'equations diff\'erentielles portant sur les applications inversibles de
$\PP^1$ dans $\PP^1$.
Ces anneaux sont munis d'une d\'erivation $D:\O_{J_k^*(\PP^1)}\rightarrow \O_{J_{k+1}^*(\PP^1)}$ d\'efinie en coordonn\'ees 
locales par :
$$
D(E)=\frac{\partial E}{\partial x}+ \frac{\partial E}{\partial y}y_1+\ldots +\frac{\partial E}{\partial y_k}y_{k+1}.
$$
Un syst\`eme d'\'equations diff\'erentielles sur les applications inversibles de $\PP^1$ dans $\PP^1$ est un faisceau 
d'id\'eaux $\I$ de $\O_{J^*(\PP^1)}$ coh\'erent (\textit{i.e.} les $\I_k = \I \cap \O_{J^*_k(\PP^1)}$ sont coh\'erents)
 diff\'erentiels (\textit{i.e.} stable par D) et r\'eduits.
 
L'espace $J_k^*(\PP^1)$ est muni d'une structure de groupo\"ide donn\'ee par : 
\begin{itemize}
\item les deux projections $s$, $t$, sur $\PP^1$,
\item une composition $c:J_k^*(\PP^1)\times_{\PP^1}J_k^*(\PP^1) \rightarrow J_k^*(\PP^1)$ d\'efinie sur les couples de jets
 $(h,g)$ tels que $t(h)=s(g)$ par les formules habituelles : 
$$
c((x,y,y_1,\ldots),(y,z,z_1,\ldots))=(x,z,z_1y_1,\ldots),
$$
\item une identit\'e $e:\PP^1 \rightarrow J_k^*(\PP^1)$ donn\'ee par $e(x)=(x,x,1,0\ldots)$,
\item une inversion $i:J_k^*(\PP^1) \rightarrow J_k^*(\PP^1)$ d\'efinie par $i(x,y,y_1\ldots)= (y,x,y_1^{-1}\ldots)$,
\end{itemize}
le tout v\'erifiant certains diagrammes commutatifs que nous ne rappellerons pas (\cite{McK}). Toutes ces fl\`eches sont 
compatibles aux faisceaux d'anneaux construits pr\'e\-c\'e\-dem\-ment dans le sens o\`u elles induisent des fl\`eches 
$s^*,t^*,e^*, i^*$ et $c^*$ entres les anneaux $\O_{\PP^1}$, $\O_{J_k^*(\PP^1)}$ et $\O_{J_k^*(\PP^1)} \otimes_{\O_{\PP^1}} 
\O_{J_k^*(\PP^1)}$ compatibles aux injections $\O_{J_k^*(\PP^1)} \subset \O_{J_{k+1}^*(\PP^1)}$.

%%%%%%%%%%%%%%%%%%%%%%%%%%%%%%%%%%%%%%%%%%%%%%%%%%%%%%%%%%%%%%%%%%%%%%%%
\subsection{$\D$-groupo\"ide de Lie}
\text{   \\}

Les d\'efinitions suivantes sont issues de \cite{malg-groupoidedegalois}.
\begin{definition}
    \label{defnaive}
Un groupo\"ide d'ordre $k$ sur $\PP^1$ est donn\'e par un faisceau d'id\'eaux coh\'erent $\mathcal{I}_k$ de 
$\mathcal{O}_{J_k^*(\PP^1)}$ tel que : 
  \begin{enumerate}
      \item $\mathcal{I}_k\subset Ker(e^*)$,
      \item $i^*\mathcal{I}_k \subset \mathcal{I}_k $,
      \item $c^*\mathcal{I}_k \subset \mathcal{I}_k\otimes_{\mathcal{O}_{\mathbb{P}^1}}1 + 1\otimes_{\mathcal{O}_{\mathbb{P}^1}}\mathcal{I}_k$ (la somme \'etant prise comme somme d'id\'eaux).
  \end{enumerate}
\end{definition}
Cette d\'efinition est naturelle mais en pratique trop restrictive pour l'utilisation que nous avons en vue. 
Une equation diff\'erentielle sur $\PP^1$ admettant des singularit\'es, l'inclusion (3) n'est pas forcement v\'erifi\'ee
au voisinage de ces points. Il faut alors utiliser la d\'efinition plus souple suivante.
\begin{definition}
    \label{defDGL}
Un $\D$-groupo\"ide de Lie sur $\mathbb{P}^1$ est donn\'e par un faisceau d'id\'eaux $\mathcal{I}$ de $\mathcal{O}_{J^*(\PP^1)}$ tel que 
  \begin{itemize}
      \item $\mathcal{I}_k=\mathcal{I}\cap \mathcal{O}_{J_k^*(\mathbb{P}^1\rightarrow\mathbb{P}^1)}$ soient coh\'erents,
      \item $\mathcal{I}$ soit stable par d\'erivation,
      \item il existe un entier $k$ et un ensemble analytique ferm\'e $Z$ dans $\PP^1$ tels que \\
(i) pour tout $\ell\geq k$, $\mathcal{I}_\ell$ v\'erifie \emph{(1)} et \emph{(2)} de la d\'efinition \ref{defnaive},\\
(ii) sur tout voisinage de $(x,y,z)\in (\PP^1-Z)\times (\PP^1-Z)\times (\PP^1-Z)$, on a \emph{(3)}.
  \end{itemize}
\end{definition}

L'exemple le plus simple est celui donn\'e par l'\'equation $S(y) = 2\frac{y_3}{y_1}-3 \left(\frac{y_2}{y_1}\right)^2= 0$. 
Les formules classiques sur la d\'eriv\'ee schwartzienne :
$$
S(c((x,y,\ldots),(y,z,\ldots))) = S(y,z,\ldots)(y_1)^2 + S(x,y,\ldots)
$$
donne les inclusions voulues. Ici l'ensemble $Z$ est vide. N\'eanmoins le th\'eor\`eme \ref{listelocaledesDGL} ci-dessous 
donne des exemples o\`u cet ensemble est non vide. Lorsque $\I$ ne contient pas d'\'equation d'ordre z\'ero, 
le $\D$-groupo\"ide de Lie est dit transitif. 
%%%%%%%%%%%%%%%%%%%%%%%%%%%%%%%%%%%%%%%%%%%%%%%%%%%%%%%%%%%%%%%%%%%%%%%%%
\subsection{$\D$-enveloppe d'une application rationnelle}
\text{   \\}

\begin{definition}
    \label{defDenv}
Soit $R:\PP^1 \to \PP^1$ une application rationnelle. Sa $\D$-enveloppe est le plus petit des $\D$-groupo\"ides de Lie 
dont $R$ est solution. Lorsque l'id\'eal de ce $\D$-groupo\"ide n'est pas $(0)$, nous dirons que la $\D$-enveloppe est 
non triviale.
\end{definition} 

L'existence d'un plus petit $\D$-groupo\"ide de Lie parmi une famille de $\D$-groupo\"ides de Lie est prouv\'ee en toute g\'en\'eralit\'e 
dans \cite{malg-groupoidedegalois} .\\

%%%%%%%%%%%%%%%%%%%%%%%%%%%%%%%%%%%%%%%%%%%%%%%%%%%%%%%%%%%%%%%%%%%%%%%%%%
%%%%%%%%%%%%%%%%%%%%%%%%%%%%%%%%%%%%%%%%%%%%%%%%%%%%%%%%%%%%%%%%%%%%%%%%%%
%%%%%%                 LISTE DES GROUPOIDES SUR PP^1                %%%%%%
%%%%%%%%%%%%%%%%%%%%%%%%%%%%%%%%%%%%%%%%%%%%%%%%%%%%%%%%%%%%%%%%%%%%%%%%%%
%%%%%%%%%%%%%%%%%%%%%%%%%%%%%%%%%%%%%%%%%%%%%%%%%%%%%%%%%%%%%%%%%%%%%%%%%%

%%%%%%%%%%%%%%%%%%%%%%%%%%%%%%%%%%%%%%%%%%%%%%%%%%%%%%%%%%%%%%%%%%%%%%%%%%
\section{Les $\D$-groupo\"ides de Lie sur $\PP^1$}
%%%%%%%%%%%%%%%%%%%%%%%%%%%%%%%%%%%%%%%%%%%%%%%%%%%%%%%%%%%%%%%%%%%%%%%%%%

Le th\'eor\`eme suivant (\cite{guy2}) donne la forme des \'equations engendrant l'id\'eal d'un $\D$-groupo\"ide de Lie 
au-dessus d'un disque $\Delta$. Pour un id\'eal $\I $ de $\O_{J^*(\PP^1)}$, nous noterons encore $\I$ l'id\'eal qu'il 
engendre dans le faisceau des \'equations diff\'erentielles m\'eromorphes ``en $x$ et en $y$'' c'est-\`a-dire dans 
$\M_{\PP^1} \otimes_{(\O_{\PP^1},s^*)} \O_{J^*(\PP^1)} \otimes_{(\O_{\PP^1},t^*)} \M_{\PP^1}$.  

\begin{theorem}
    \label{listelocaledesDGL}
Soit $\I$ l'id\'eal d'un \dgl transitif au-dessus d'un disque $\Delta$. Il est diff\'e\-rentia\-blement engendr\'e par 
une seule \'equation m\'eromorphe d'ordre inf\'erieur ou \'egal \`a trois d'une des quatres formes suivantes:  
  \begin{enumerate}
      \item[(1)] $\eta(y)(y_1)^n-\eta(x)=0$ avec $n$ entier et $\eta$ m\'eromorphe sur $\Delta$, 
      \item[(2)] $\mu(y)y_1 + \frac{y_2}{y_1}-\mu(x)=0 $ avec $\mu$ m\'eromorphe,
      \item[(3)] $\nu(y)(y_1)^2 +2\frac{y_3}{y_1}-3 \left( \frac{y_2}{y_1}\right)^2-\nu(x)=0$ avec $\nu$ m\'eromorphe, 
      \item[($\infty$)] $0=0$ 
  \end{enumerate}
Les \dgls correspondants seront respectivement not\'e $G_1^n(\eta)$, $G_2(\mu)$, $G_3(\nu)$ et $G_\infty$.
\end{theorem}
Lorsqu'on effectue un changement de coordonn\'ee, ces \'equations sont soumis aux transformations suivantes.
\begin{proposition}
    \label{transfodejauge}
Soient $\varphi : \Delta_1 \to \Delta_2 $ un diff\'eomorphisme et $x_i$ une coordonn\'ee sur $\Delta_i$. L'application 
$\varphi$ donne naturellement une application $\varphi_* : J_k^*(\Delta_1) \to J_k^*(\Delta_2)$. L'image r\'eciproque par 
$\varphi$ d'un \dgl sur $\Delta_2$ est un \dgl sur $\Delta_2$. Il est donn\'e par l'une des \'equations suivantes : 
$$
\begin{array}{l}                                           
      \varphi^*  G_0(h) = G_0(h\circ\varphi)     \\
      \varphi^*  G_1^n(\eta) = G_1^n(\eta\circ\varphi (\varphi')^n)     \\
      \varphi^*  G_2(\mu) = G_2(\mu\circ\varphi \varphi'+\frac{\varphi''}{\varphi'})      \\
      \varphi^*  G_3(\nu) = G_3(\nu\circ\varphi(\varphi')^2 + S(\varphi))                \\ 
\end{array}
$$
o\`u $S(\varphi)$ est la d\'eriv\'ee Schwartzienne de $\varphi$ par rapport \`a la coordonn\'ee $x_1$.
\end{proposition}

Nous allons nous int\'eresser aux \'equations d'un \dgl sur $\PP^1$. Le r\'esultat suivant est un corollaire du th\'eor\`eme
\ref{listelocaledes DGL}.

\begin{proposition}
Soit $\I$ l'id\'eal d'un \dgl transitif au dessus de $\PP^1$. Il est diff\'erentiablement engendr\'e par une des 
\'equations suivantes : 

\begin{enumerate}
  \item[(1)] $\eta(y)(y_1)^n-\eta(x)=0$ avec $n$ entier,

  \item[(2)] $\mu(y)y_1 + \frac{y_2}{y_1}-\mu(x)=0 $,
  \item[(3)] $\nu(y)(y_1)^2 +2\frac{y_3}{y_1}-3 \left( \frac{y_2}{y_1}\right)^2-\nu(x)=0$, 
  \item[($\infty$)] $0=0$ .

\end{enumerate}
 avec $\eta$, $\mu$ et $\nu$ rationnelles.
\end{proposition}
\begin{preuve}
Montrons cette proposition dans le cas (3), les autres cas se montre de la m\^eme mani\`ere. Sur chaque ``disque'' 
 $\mathcal{U}=\PP^1 - \{\infty\}$ (de coordonn\'ee $x$) et $\mathcal{V}=\PP^1 - \{0\}$ (de coordonn\'ee 
 $\overline{x}=\frac{1}{x}$), l'id\'eal d'un \dgl d'ordre trois est engendr\'e par une \'equation de type (3) : $G_3(\nu)$
 sur $\mathcal{U}$ et $G_3(\overline {\nu})$ sur $\mathcal{V}$. D'apr\`es la proposition \ref{transfodejauge}, ces deux 
 \'equations sont reli\'ees par l'identit\'e : $\overline{\nu}(\frac{1}{x})\frac{1}{x^4}=\nu (x)$. Ceci prouve la rationnalit\'e de $\nu$. 
    
Remarquons qu'on peut \'etablir une preuve ind\'ependante du th\'eor\`eme \ref{listelocaledesDGL} \`a partir du fait qu'un \dgl transitif est le groupo\"ide d'invariance d'une structure g\'eom\'etrique et que 
celle-ci soit complètement d\'etermin\'ee par le groupe d'isotropie d'un point. Dans le cas non-singulier, nous renvoyons 
le lecteur \`a \cite{koba} et pour l'adaptation au cas rationnel \`a \cite{mathese}.  
\end{preuve}

%%%%%%%%%%%%%%%%%%%%%%%%%%%%%%%%%%%%%%%%%%%%%%%%%%%%%%%%%%%%%%%%%%%%%%%%%%%%%%%%%%%%%
%%%%%%%%%%%%%%%%%%%%%%%%%%%%%%%%%%%%%%%%%%%%%%%%%%%%%%%%%%%%%%%%%%%%%%%%%%%%%%%%%%%%%
%%%%%%                             LA PREUVE                                   %%%%%%
%%%%%%%%%%%%%%%%%%%%%%%%%%%%%%%%%%%%%%%%%%%%%%%%%%%%%%%%%%%%%%%%%%%%%%%%%%%%%%%%%%%%%
%%%%%%%%%%%%%%%%%%%%%%%%%%%%%%%%%%%%%%%%%%%%%%%%%%%%%%%%%%%%%%%%%%%%%%%%%%%%%%%%%%%%%

%%%%%%%%%%%%%%%%%%%%%%%%%%%%%%%%%%%%%%%%%%%%%%%%%%%%%%%%%%%%%%%%%%%%%%%%%%%%%%%%%%%%%
\section{Preuve du th\'eor\`eme principal}
%%%%%%%%%%%%%%%%%%%%%%%%%%%%%%%%%%%%%%%%%%%%%%%%%%%%%%%%%%%%%%%%%%%%%%%%%%%%%%%%%%%%%

Remarquons d'abord que la $\mathcal{D}$-enveloppe d'une homographie est non triviale. En effet, dans ce cas, $R$ est 
solution de $G_3(0)$ d'\'equation $S(y)=0$.
Nous supposerons donc que $R$ n'est pas une homographie.\\

Dans un premier temps, nous allons supposer que $R$ est solution d'une \'equation diff\'erentielle $G_2(\mu)$ et \'etablir 
le lemme suivant.

\begin{lemma}
Une transformation rationnelle $R$ de $\PP^1$ est solution d'un \dgl $G_2(\mu)$ si et seulement si $R$ provient du passage au 
quotient par un sous-groupe d'application affine d'une dilatation de $\CC$.
\end{lemma}

\begin{preuve} 
Quitte \`a remplacer $R$ par $R^{\circ n}$, nous pouvons trouver un point $p \in \mathbb{P}^1$ en dehors des p\^oles de 
 $\mu$, fixe et r\'epulsif.
Le th\'eor\`eme de K{\oe}nigs nous permet alors de construire une lin\'earisante locale holomorphe que nous \'etendons \`a $\CC$ :
$$
\Psi : \mathbb{C} \rightarrow \mathbb{P}^1-E \text{ v\'erifiant } \Psi(\lambda z)=R \circ \Psi \text{ et } \Psi'(0)=1.
$$ 
D'apr\`es la proposition \ref{transfodejauge}, l'image r\'eciproque par $\Psi$ du  \dgl $G_2(\mu)$ est le \dgl au-dessus
 de $\CC$ : $G_2(\overline{\mu})$ o\`u$$
\overline{\mu}=(\mu \circ \Psi) \Psi' + \frac{\Psi''}{\Psi'}.
$$ 
La fraction $R$ \'etant solution de $G_2(\mu)$, l'homoth\'etie $z\to\lambda z$ est solution de  $G_2(\overline{\mu})$ d'o\`u 
$$\overline{\mu}(\lambda z) \lambda = \overline{\mu}(z).$$
Comme $|\lambda|>1$, l'\'egalit\'e ci-dessus implique l'existence d'une constante $c$ telle que 
$\overline{\mu}(z)=\frac{c}{z}$. Ayant suppos\'e $\mu(p)$ fini et $\Psi'(0)=1 $, $\overline{\mu}(0)$ doit \^etre fini ce qui 
force $\overline{\mu}$ \`a \^etre nul. L'image r\'eciproque de $G_2(\mu)$ est donc le \dgl $G_2(0)$, d'\'equation 
$\frac{y_2}{y_1}=0$,
 dont les solutions sont les applications affines de $\CC$.
Consid\'erons deux point $p$ et $q$ de $\CC$ qui ne sont pas des points critiques de $\Psi$ tels que $\Psi(p)=\Psi(q)$. 
Il existe une application locale $\gamma : (\CC,p) \to (\CC,q)$ v\'erifiant $\Psi \circ \gamma = \Psi$. 
Par construction de $\overline{\mu}$, cette application $\gamma$ laisse le \dgl $G_2(\overline{\mu})$ invariant.
 La proposition \ref{transfodejauge} $\gamma$ est solution de $G_2(\overline{\mu})$ 
 donc est une application affine.

Le groupe $G=\{ \gamma \ | \ \Psi\circ\gamma=\Psi \}$ est un groupe d'applications affines qui agit transitivement sur 
les fibres de $\Psi$. Ces fibres \'etant discr\`etes, ce groupe est discret.

Nous avons une condition n\'ecessaire pour que la $\D$-enveloppe de $R$ soit de la forme $G_2(\mu)$.\\

Montrons que cette condition est aussi suffisante. 
Consid\'ereons $\Psi : \CC \to \PP1$ une lin\'earisante de K\oe nigs de $R$ et $G_2(0)$ le groupo\"ide des transformations affines
de $\CC$.
Localement, on peut ramener par une d\'etermination de $\Psi^{-1}$ le \dgl $G_2(0)$ en un \dgl au-dessus d'un ouvert de $\PP^1$.
 D'apr\`es
la proposition \ref{transfodejauge}, ce \dgl est $G_2(\frac{(\Psi^{-1})''}{(\Psi^{-1})'})$ et $R$ en est une 
solution (lorsque cela \`a un sens). \\
Il nous suffit de v\'erifier que l'\'equation de ce groupo\"ide est en fait rationnelle, c'est-\`a-dire, pour tous les 
groupes $G$ que l'on peut rencontrer, v\'erifier que $\frac{(\Psi^{-1})''}{(\Psi^{-1})'}$ est rationnelle.

 Consid\'erons le r\'eseau $\Lambda=\{b \ | \ (z \mapsto z+b) \in G\}$. 
On obtient une factorisation de $\Psi$:
$$
\mathbb{C} \rightarrow \mathbb{C}/\Lambda \rightarrow \mathbb{P}^1
$$ 
La derni\`ere application induite par $\Psi$ est une fonction m\'eromorphe sur $\mathbb{C}/\Lambda$ invariante sous l'action 
de $G/\Lambda$. 
En utilisant la liste des sous-groupes discrets d'applications affines, on obtient la liste suivante d'applications :   

\begin{enumerate}

\item $\Lambda=\{0\}$ et $G$ est un groupe fini de rotation : \\ 
$\Psi(z)=z^k$ , $R$ est une homoth\'etie et $\mu=0$. 

\item $\Lambda=\mathbb{Z}$ et $G=\Lambda$ : \\ 
$\Psi(z)= \exp(2i\pi z)$, $R$ est un mon\^ome et $\mu(z)=\frac{-1}{z}$. 

\item $\Lambda=\mathbb{Z}$ et $G/\Lambda=\{+1,-1\}$ : \\ 
$\Psi(z)= \cos(2i\pi z)$, $R$ est un polyn\^ome Tch\'ebitchev et $\mu(z)=\frac{-z}{z^2-4}$. 

\item $\Lambda=\mathbb{Z}+\mathbb{Z}\tau$ ($\Im \tau > 0$) et $G/\Lambda=\{+1,-1\}$ : \\ 
$\Psi(z)=\wp(z)$ ($\wp$ est la fonction de Weierstrass associ\'ee au r\'eseau $\Lambda$ 
et solution de $(\wp')^2=4\wp^3+g_2\wp+g_3$ pour des constantes $g_2$ et $g_3$ d\'efinies par $\Lambda$ ), $R$ est un exemple de Latt\`es et $\mu(z)=-\frac{6z^2+g_2/2}{4z^3+g_2z+g_3}$.

\item $\Lambda=\mathbb{Z}+\mathbb{Z}i$ et $G/\Lambda=\{+1,i,-1,-i\}$ : \\ 
dans ce cas $g_3=0$, $\Psi(z)=\wp(z)^2$ et $\mu(z)=-\frac{1}{4z}-\frac{6z+g_2/2}{8z^2+2g_2z}$.

\item $\Lambda=\mathbb{Z}+\mathbb{Z}j$ et $G/\Lambda=\{+1,j,j^2\}$ : \\ 
dans ce cas $g_2=0$, $\Psi(z)=\wp'(z)$ et $\mu(z)=-\frac{2}{3}\frac{z}{z^2-g_3)}$.

\item $\Lambda=\mathbb{Z}+\mathbb{Z}j$ et $G/\Lambda=\{+1,j,j^2,-j,-j^2,-1\}$ : \\ 
dans ce cas $g_2=0$, $\Psi(z)=\wp(z)^3$ et $\mu(z)=-\frac{2}{9z}-\frac{1}{z+g_3/2}$.

\end{enumerate} 

\text{  \\  }

\end{preuve}

\begin{lemma}
Aucune transformation rationnelle de $\PP^1$ n'a de $\mathcal{D}$-enveloppe de type (3)
\end{lemma}

\begin{preuve}
Supposons que $R$ soit solution d'un \dgl $G_3(\nu)$. En proc\'edant comme pr\'ec\'edemment, on obtient une lin\'earisante 
locale $\Psi$ par le th\'eor\`eme de K\oe nigs que l'on prolonge en un rev\^etement $\CC \to \PP^1$. 
L'homoth\'etie $z \to \lambda z$ \'etant solution de l'image r\'eciproque de $G_3(\nu)$ par $\Psi$ : $G_3(\overline{\nu})$. 
On a l'\'egalit\'e :
$$
\overline{\nu}(\lambda z) \lambda^2 = \overline{\nu}(z),
$$
d'o\`u on d\'eduit qu'il existe une constante $c$ telle que $\overline{\nu}= \frac{c}{z^2}$. Comme $\overline{\nu}(0)$ est 
fini, on a $\overline{\nu}=0$. L'image r\'eciproque de $G_3(\nu)$ est le groupo\"ide des homographies. 
Soit $\gamma$ un germe d'application v\'erifiant $\Psi \circ \gamma= \Psi$. On \'etablit comme pr\'ec\'edemment 
que $S(\gamma)=0$ et donc que $\gamma$ est une homographie. Il existe donc un sous-groupe $G$ des homographie qui agit transitivement 
sur les fibres de $\Psi$.
Ces fibres de $\Psi$ \'etant des sous ensembles discrets de $\CC$, $G$ ne peut \^etre compos\'e que d'applications affines.
On obtient la m\^eme liste d'application que pour les solutions de \dgl de type (2). On en d\'eduit qu'aucune application rationnelle n'a de $\D$-enveloppe de type (3).    

\end{preuve}

\section{remarque}

1 - Une caract\'erisation \textit{a priori} partant sur le $\mathcal{D}$-enveloppe de $R$ de l'existence d'un 
commutant non trivial est n\'ecessaire pour aborder l'\'etude des applications rationnelles en dimension sup\'erieure.
Une condition n\'ecessaire semble \^etre que la $\mathcal{D}$-enveloppe soit localement l'action d'un groupe sur un espace 
homog\`ene. Cette condition n'est pas suffisante au vue des travaux de \cite{D-S}. Le groupe de Lie semble devoir \^etre un groupe
de transformation affine (voir aussi \cite{Veselov}).

2 - Sur une courbe alg\'ebrique, l'\'etude de la $\mathcal{D}$-enveloppe d'une correspondance reste \`a faire (\cite{CU}). Il est
trivial de v\'erifier que les correspondances modulaires ont une $\mathcal{D}$-enveloppe non triviale. En utilisant les r\'esultats
de \cite{margulis}, on v\'erifie que si la $\mathcal{D}$-enveloppe d'une correspondance est non triviale et non singuli\`ere 
(\textit{i.e.} $G_3(\nu)$ avec $\nu$ holomorphe) alors
la correspondance est modulaire. Le cas des correspondances ayant une $\mathcal{D}$-enveloppe singuli\`ere reste obscur.

3 - L'\'equation de la lin\'earisante de K\oe nigs est une \'equation aux $q$-diff\'erences non lin\'eaire et autonome. Une 
\'etude de la $\mathcal{D}$-enveloppe des syst\`emes dynamiques provenant d'\'equations aux diff\'erences ou aux $q$-diff\'erences
permettrait de relier l'\textit{``int\'egrabilit\'e''} d'une transformation rationnel \`a l'\textit{``int\'egrabilit\'e''} de 
ses lin\'earisantes de K\oe nigs (ou nor\-malisantes de Poincar\'e).

%%%%%%%%%%%%%%%%%%%%%%%%%%%%%%%%%%%%%%%%%%%%%%%%%%%%%%%%%%%%%%%%%%%%%%%%%%%%%
%%%%%%%%%%%%%%%%%%%%%%%%%%%%%%%%%%%%%%%%%%%%%%%%%%%%%%%%%%%%%%%%%%%%%%%%%%%%%
%%%%%                     EN DIMENSION SUPERIEUR                       %%%%%%
%%%%%%%%%%%%%%%%%%%%%%%%%%%%%%%%%%%%%%%%%%%%%%%%%%%%%%%%%%%%%%%%%%%%%%%%%%%%%
%%%%%%%%%%%%%%%%%%%%%%%%%%%%%%%%%%%%%%%%%%%%%%%%%%%%%%%%%%%%%%%%%%%%%%%%%%%%%

%%%%%%%%%%%%%%%%%%%%%%%%%%%%%%%%%%%%%%%%%%%%%%%%%%%%%%%%%%%%%%%%%%%%%%%%%%%%%  
%\section{En dimension sup\'erieur}
%%%%%%%%%%%%%%%%%%%%%%%%%%%%%%%%%%%%%%%%%%%%%%%%%%%%%%%%%%%%%%%%%%%%%%%%%%%%%

\end{document}